\documentclass[11pt]{article}

\usepackage{amsfonts,amssymb,amsmath,graphicx,graphics}
\usepackage[english]{babel}
\makeatletter
\@addtoreset{equation}{section}
\makeatother

\newcounter{enunciato}[section]

\newtheorem{ittheorem}{Theorem}
\newtheorem{itlemma}{Lemma}
\newtheorem{itproposition}{Proposition}

\newenvironment{theorem}{\addtocounter{enunciato}{1}
\begin{ittheorem}}{\end{ittheorem}}

\newenvironment{proposition}{\addtocounter{enunciato}{1}
\begin{itproposition}}{\end{itproposition}}

\def \cW {{\mathcal W}}
\def \cD {{\mathcal D}}
\def \cL {{\mathcal L}}
\def \cI {{\mathcal I}}
\def \cR {{\mathcal R}}
\def \cA {{\mathcal A}}
\def \CONE {{\hbox{\footnotesize\rm CONE}}}

\begin{document}

\title{A mathematical model for a copolymer in an emulsion}

\author{\renewcommand{\thefootnote}{\arabic{footnote}}
F.\ den Hollander
\footnotemark[1]\,\,\,\footnotemark[2]
\\
\renewcommand{\thefootnote}{\arabic{footnote}}
N.\ P\'etr\'elis
\footnotemark[2]
}

\footnotetext[1]
{Mathematical Institute, Leiden University, P.O.\ Box 9512,
2300 RA Leiden, The Netherlands}\,

\footnotetext[2]
{EURANDOM, P.O.\ Box 513, 5600 MB Eindhoven, The Netherlands}

\maketitle

\begin{abstract}
In this paper we review some recent results, obtained jointly with Stu
Whittington, for a mathematical model describing a copolymer in an emulsion.
The copolymer consists of hydrophobic and hydrophilic monomers, concatenated
randomly with equal density. The emulsion consists of large blocks of oil
and water, arranged in a percolation-type fashion. To make the model
mathematically tractable, the copolymer is allowed to enter and exit a
neighboring pair of blocks only at diagonally opposite corners. The energy
of the copolymer in the emulsion is minus $\alpha$ times the number of
hydrophobic monomers in oil minus $\beta$ times the number of hydrophilic
monomers in water. Without loss of generality we may assume that the
interaction parameters are restricted to the cone $\{(\alpha,\beta)\in
\mathbb{R}^2\colon\,|\beta|\leq\alpha\}$.

We show that the phase diagram has two regimes: (1) in the supercritical
regime where the oil blocks percolate, there is a single critical curve in
the cone separating a localized and a delocalized phase; (2) in the subcritical
regime where the oil blocks do not percolate, there are three critical
curves in the cone separating two localized phases and two delocalized
phases, and meeting at two tricritical points. The different phases are
characterized by different behavior of the copolymer inside the four
neighboring pairs of blocks.

\vskip 0.5truecm
\noindent
\emph{AMS} 2000 \emph{subject classifications.} 60F10, 60K37, 82B27.\\
\emph{Key words and phrases.} Random copolymer, random emulsion, localization,
delocalization, phase transition, percolation.\\
{\it Acknowledgment.} NP is supported by a postdoctoral fellowship from the
Netherlands Organization for Scientific Research (grant 613.000.438).

\vspace{1cm}
\noindent
* Invited paper to appear in a special volume of the Journal of Mathematical
Chemistry on the occasion of the 65th birthdays of Ray Kapral and Stu Whittington
from the Department of Chemistry at the University of Toronto.
\end{abstract}

\newpage


\section{Introduction}
\label{S1}

\begin{figure}
\begin{center}
\includegraphics[scale = 0.3]{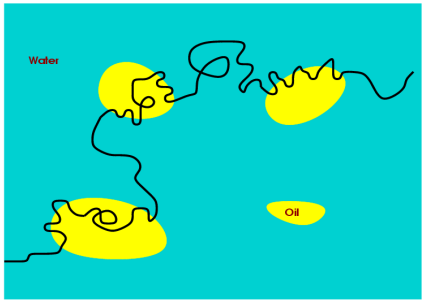}
\end{center}
\caption{An undirected copolymer in an emulsion.}
\label{fig-PolInRa}
\end{figure}

The physical situation we want to model is that of a copolymer in an emulsion
(see Fig.~\ref{fig-PolInRa}). The random interface model described below was
introduced in den Hollander and Whittington~\cite{dHoWh06}, where the qualitative
properties of the phase diagram were obtained. Finer details of the phase diagram
are derived in den Hollander and P\'etr\'elis~\cite{dHoPe07a}, \cite{dHoPe07b}.

\subsection{The model}
\label{S1.1}

\begin{figure}
\begin{center}
\includegraphics[scale = 0.4]{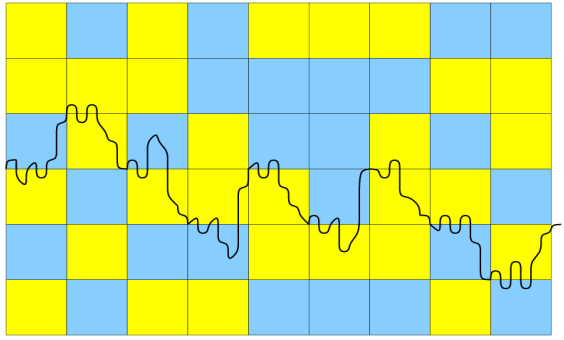}
\end{center}
\caption{A directed self-avoiding path crossing blocks of oil and water diagonally.
The light-shaded blocks are oil, the dark-shaded blocks are water. Each block is
$L_n$ lattice spacings wide in both directions. The path carries hydrophobic and
hydrophilic monomers on the lattice scale, which are not indicated.}
\label{fig-copolemulblock}
\end{figure}

Each positive integer is randomly labelled $A$ or $B$, with probability
$\frac{1}{2}$ each, independently for different integers. The resulting
labelling is denoted by
\begin{equation}
\label{bondlabel}
\omega = \{\omega_i \colon\, i \in \mathbb{N}\} \in \{A,B\}^\mathbb{N}
\end{equation}
and represents the \emph{randomness of the copolymer}, with $A$ indicating
a hydrophobic monomer and $B$ a hydrophilic monomer. Fix $p \in (0,1)$ and
$L_n \in \mathbb{N}$. Partition $\mathbb{R}^2$ into square blocks of size
$L_n$:
\begin{equation}
\label{blocks}
\mathbb{R}^2 = \bigcup_{x \in \mathbb{Z}^2} \Lambda_{L_n}(x), \qquad
\Lambda_{L_n}(x) = xL_n + (0,L_n]^2.
\end{equation}
Each block is randomly labelled $A$ or $B$, with probability $p$, respectively,
$1-p$, independently for different blocks. The resulting labelling is denoted by
\begin{equation}
\label{blocklabel}
\Omega = \{\Omega(x) \colon\, x \in \mathbb{Z}^2\} \in \{A,B\}^{\mathbb{Z}^2}
\end{equation}
and represents the \emph{randomness of the emulsion}, with $A$ indicating an
oil block and $B$ a water block (see Fig.~\ref{fig-copolemulblock}).

Let $\cW_n$ be the set of $n$-step \emph{directed self-avoiding paths} starting at
the origin and being allowed to move \emph{upwards, downwards and to the right}. Let
\begin{itemize}
\item[$\bullet$]
$\cW_{n,L_n}$ is the subset of $\cW_n$ consisting of those paths that enter blocks at
a corner, exit blocks at one of the two corners \emph{diagonally opposite} the one
where it entered, and in between \emph{stay confined} to the two blocks that are
seen upon entering.
\end{itemize}
In other words, after the path reaches a site $xL_n$ for some $x\in\mathbb{Z}^2$,
it must make a step to the right, it must subsequently stay confined to the pair
of blocks labelled $x$ and $x+(0,-1)$, and it must exit this pair of blocks either
at site $xL_n+(L_n,L_n)$ or at site $xL_n+(L_n,-L_n)$ (see Fig.~\ref{fig-copolemulblock}).
This restriction -- which is put in to make the model mathematically tractable -- is
unphysical. Nonetheless, the model has physically very relevant behavior.

Given $\omega,\Omega$ and $n$, with each path $\pi \in \cW_{n,L_n}$ we associate
an energy given by the Hamiltonian
\begin{equation}
\label{Hamiltonian}
H_{n,L_n}^{\omega,\Omega}(\pi)
= - \sum_{i=1}^n
\left(\alpha\,1\left\{\omega_i=\Omega^{L_n}_{(\pi_{i-1},\pi_i)}=A\right\}
+ \beta\,1\left\{\omega_i=\Omega^{L_n}_{(\pi_{i-1},\pi_i)}=B\right\}\right),
\end{equation}
where $\alpha,\beta\in\mathbb{R}$ and $\Omega^{L_n}_{(\pi_{i-1},\pi_i)}$ denotes
the label of the block that the edge $(\pi_{i-1},\pi_i)$ lies in. What this Hamiltonian
does is count the number of $AA$-matches and $BB$-matches and assign them energy
$-\alpha$ and $-\beta$, respectively. Note that the interaction is assigned to
edges rather than to vertices, i.e., we identify the monomers with the steps of
the path. We will see shortly that without loss of generality we may restrict the
interaction parameters to the cone
\begin{equation}
\label{defcone}
\CONE = \{(\alpha,\beta)\in\mathbb{R}^2\colon\,\alpha\geq |\beta|\}.
\end{equation}

\subsection{The free energy}
\label{S1.2}

Given $\omega,\Omega$ and $n$, we define the \emph{quenched free energy per step}
as
\begin{equation}
\label{fedef}
\begin{aligned}
f_{n,L_n}^{\omega,\Omega}
&= \frac{1}{n} \log Z_{n,L_n}^{\omega,\Omega},\\
Z_{n,L_n}^{\omega,\Omega}
&= \sum\limits_{\pi\in\cW_{n,L_n}} \exp\left[-H_{n,L_n}^{\omega,\Omega}(\pi)\right].
\end{aligned}
\end{equation}
We are interested in the limit $n \to \infty$ subject to the restriction
\begin{equation}
\label{Ln}
L_n \to \infty \qquad \mbox{ and } \qquad \frac{1}{n}L_n \to 0.
\end{equation}
This is a \emph{coarse-graining} limit where the path spends a long time in each
single block yet visits many blocks. In this limit, there is a separation between
a \emph{polymer scale} and an \emph{emulsion scale} (see Fig.~\ref{fig-copolemulblock}).

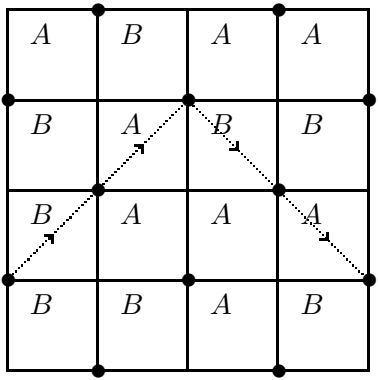
\begin{figure}
\begin{center}
\setlength{\unitlength}{0.3cm}
\begin{picture}(10,10)(4,0)
{\thicklines
\qbezier(0,16)(8,16)(16,16)
\qbezier(16,0)(16,8)(16,16)
\qbezier(0,0)(8,0)(16,0)
\qbezier(0,0)(0,8)(0,16)
\qbezier(4,0)(4,8)(4,16)
\qbezier(8,0)(8,8)(8,16)
\qbezier(12,0)(12,8)(12,16)
\qbezier(0,4)(8,4)(16,4)
\qbezier(0,8)(8,8)(16,8)
\qbezier(0,12)(8,12)(16,12)
}
\qbezier[30](0,4)(2,6)(4,8)
\qbezier[30](4,8)(6,10)(8,12)
\qbezier[30](8,12)(10,10)(12,8)
\qbezier[30](12,8)(14,6)(16,4)
{\thicklines
\qbezier(1.7,6)(1.85,6)(2,6)
\qbezier(2,5.7)(2,5.85)(2,6)
\qbezier(5.7,10)(5.85,10)(6,10)
\qbezier(6,9.7)(6,9.85)(6,10)
\qbezier(9.9,9.8)(10.05,9.8)(10.2,9.8)
\qbezier(10.2,10.1)(10.2,9.95)(10.2,9.8)
\qbezier(13.9,5.8)(14.05,5.8)(14.2,5.8)
\qbezier(14.2,6.1)(14.2,5.95)(14.2,5.8)
}
\put(1,14.5){$A$}
\put(5,14.5){$B$}
\put(9,14.5){$A$}
\put(13,14.5){$A$}
\put(1,10.5){$B$}
\put(5,10.5){$A$}
\put(9,10.5){$B$}
\put(13,10.5){$B$}
\put(1,6.5){$B$}
\put(5,6.5){$A$}
\put(9,6.5){$A$}
\put(13,6.5){$A$}
\put(1,2.5){$B$}
\put(5,2.5){$B$}
\put(9,2.5){$A$}
\put(13,2.5){$B$}
\put(0.05,4){\circle*{.6}}
\put(0.05,12){\circle*{.6}}
\put(4.05,0){\circle*{.6}}
\put(4.05,8){\circle*{.6}}
\put(4.05,16){\circle*{.6}}
\put(8.05,4){\circle*{.6}}
\put(8.05,12){\circle*{.6}}
\put(12.05,0){\circle*{.6}}
\put(12.05,8){\circle*{.6}}
\put(12.05,16){\circle*{.6}}
\put(16.05,4){\circle*{.6}}
\put(16.05,12){\circle*{.6}}
\end{picture}
\end{center}
\caption{The coarse-grained path sampling $\Omega$. The dashed lines with
arrows, which denote the steps in this path, indicate where the copolymer
enters and exits. In between, the copolymer stays confined to the two
neighboring blocks, as in Fig.~\ref{fig-copolemulblock}. The block pairs
visited by the coarse-grained path in this example are $BB$, $AA$, $BA$
and $AB$, where the first (second) symbol indicates the type of block that
is crossed (not crossed) diagonally.}
\label{fig-coarsesampl}
\end{figure}

The starting point of the analysis is the following \emph{variational representation}
of the free energy. Let
\begin{itemize}
\item
$\cA$ is the set of all $2 \times 2$-matrices $(a_{kl})_{kl\in\{A,B\}^2}$ whose
elements are $\geq 2$.
\item
$\cR(p)$ is the set of all $2 \times 2$ matrices $(\rho_{kl})_{kl\in\{A,B\}^2}$
whose elements are the possible \emph{limiting frequencies} at which the four types
of block pairs are visited along a \emph{coarse-grained path} (= a path on the
corners of the blocks crossing diagonally), with $k$ indicating the label of the
block that is crossed and $l$ indicating the label of the block that is not crossed
(see Fig.~\ref{fig-coarsesampl}).
\item
$(\psi_{kl}(\alpha,\beta;a_{kl}))_{kl\in\{A,B\}^2}$ is the $2 \times 2$ matrix
of free energies per step of the copolymer in a $kl$-block of size $L \times L$
when the total number of steps inside the block is $a_{kl}L$, in the limit as
$L\to\infty$.
\end{itemize}

\begin{theorem}
\label{feiden}
For all $(\alpha,\beta)\in\mathbb{R}^2$ and $p \in (0,1)$,
\begin{equation}
\label{sa}
\lim_{n\to\infty} f_{n,L_n}^{\omega,\Omega} = f = f(\alpha,\beta;p)
\end{equation}
exists $\omega,\Omega$-a.s., is finite and non-random, and is given by
\begin{equation}
\label{fevar}
f = \sup_{(a_{kl}) \in \cA}\, \sup_{(\rho_{kl}) \in \cR(p)}
\frac{\sum_{kl} \rho_{kl} a_{kl} \psi_{kl}(\alpha,\beta;a_{kl})}
{\sum_{kl} \rho_{kl} a_{kl}}.
\end{equation}
\end{theorem}

\noindent
Theorem~\ref{feiden} says that the free energy per step is obtained by keeping
track of the times spent in each of the four types of block pairs, summing the
free energies of the four types of block pairs given these times, and afterward
optimizing over these times and over the coarse-grained random walk. Note that
the latter carries no entropy, because of (\ref{Ln}). For details of the proof
we refer to den Hollander and Whittington~\cite{dHoWh06}.

It can be shown that $f(\alpha,\beta;p)$ is convex in $(\alpha,\beta)$ and
continuous in $p$, and has the symmetry properties
\begin{equation}
\label{symms}
\begin{aligned}
f(\alpha,\beta;p) &= f(\beta,\alpha;1-p),\\
f(\alpha,\beta;p) &= \textstyle{\frac12}(\alpha+\beta) + f(-\beta,-\alpha;p).
\end{aligned}
\end{equation}
\noindent
These are the reason why without loss of generality we may restrict the
parameters to the cone in (\ref{defcone}).

Theorem~\ref{feiden} shows that, in order to get the phase diagram, what we
need to do is collect the necessary information on the \emph{two key ingredients}
of (\ref{fevar}), namely, the four block pair free energies $\psi_{kl}$,
$k,l\in\{A,B\}$, and the percolation set $\cR(p)$. We will see that only very
little is needed about $\cR(p)$.

The behavior of $f$ as a function of $(\alpha,\beta)$ is different for $p
\geq p_c$ and $p < p_c$, where $p_c \approx 0.64$ is \emph{the critical
percolation density for directed bond percolation on the square lattice}.
The reason is that the coarse-grained paths, which determine the set $\cR(p)$,
sample $\Omega$ just like paths in directed bond percolation on the square
lattice rotated by 45 degrees sample the percolation configuration (see
Fig.~\ref{fig-coarsesampl}).

\subsection{Free energies per pair of blocks}
\label{S1.3}

Because $AA$-blocks and $BB$-blocks have no interface, we have for all
$(\alpha,\beta)\in\mathbb{R}^2$ and $a\geq 2$,
\begin{equation}
\label{psiAABB}
\psi_{AA}(\alpha,\beta;a) = \textstyle{\frac12}\alpha+\kappa(a,1)
\qquad \mbox{ and } \qquad
\psi_{BB}(\alpha,\beta;a) = \textstyle{\frac12}\beta+\kappa(a,1),
\end{equation}
where $\kappa(a,1)$ is the entropy per step of walks that diagonally cross
a block of size $L \times L$ in $aL$ steps, in the limit as $L\to\infty$.
There is an explicit formula for $\kappa(a,1)$, which we will not specify
here.

To compute $\psi_{AB}(\alpha,\beta;a)$ and $\psi_{BA}(\alpha,\beta;a)$ is
harder. The following \emph{variational representation} holds.

\begin{proposition}
\label{p:linkinf}
For all $(\alpha,\beta)\in\mathbb{R}^2$ and $a\geq 2$,
\begin{equation}
\label{psiinflink}
a\psi_{AB}(\alpha,\beta;a) = \sup_{ {0 \leq b \leq 1,\,c\geq b} \atop {a-c \geq 2-b} }
\left\{c\phi^{\cI}(\alpha,\beta;c/b)
+(a-c)\left[\tfrac{1}{2}\alpha+\kappa(a-c,1-b)\right]\right\},
\end{equation}
where $\phi^\cI(\alpha,\beta;c/b)$ is the free energy per step associated with
walks running along a linear interface over a distance $cL$ in $bL$ steps, and
$\kappa(a-c,1-b)$ is the entropy per step of walks that diagonally cross a block
of size $(1-b)L \times L$ in $(a-c)L$ steps, both in the limit as $L\to\infty$.
\end{proposition}

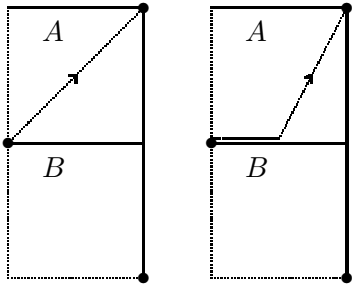
\begin{figure}
\begin{center}
\setlength{\unitlength}{0.3cm}
\begin{picture}(10,10)(2,-3)
{\thicklines
\qbezier(0,6)(3,6)(6,6)
\qbezier(6,0)(6,3)(6,6)
\qbezier(0,0)(3,0)(6,0)
\qbezier(6,-6)(6,-3)(6,0)
}
\qbezier[40](0,0)(0,3)(0,6)
\qbezier[40](0,0)(0,-3)(0,-6)
\qbezier[40](0,-6)(3,-6)(6,-6)
{\thicklines
\qbezier(2.7,3)(2.85,3)(3,3)
\qbezier(3,2.7)(3,2.85)(3,3)
\qbezier(13.15,2.9)(13.3,2.95)(13.45,3)
\qbezier(13.45,2.7)(13.45,2.85)(13.45,3)
}
\put(1.5,4.5){$A$}
\put(1.5,-1.5){$B$}
\qbezier[60](0,0)(3,3)(6,6)
{\thicklines
\qbezier(9,6)(12,6)(15,6)
\qbezier(15,0)(15,3)(15,6)
\qbezier(9,0)(12,0)(15,0)
\qbezier(15,-6)(15,-3)(15,0)
}
\qbezier[40](9,0)(9,3)(9,6)
\qbezier[40](9,0)(9,-3)(9,-6)
\qbezier[40](9,-6)(12,-6)(15,-6)
\put(10.5,4.5){$A$}
\put(10.5,-1.5){$B$}
\qbezier[50](12,0.2)(13.5,3.1)(15,6)
\qbezier[30](9,.2)(10.5,.2)(12,.2)
\put(0,0){\circle*{.5}}
\put(6,6){\circle*{.5}}
\put(6,-6){\circle*{.5}}
\put(9,0){\circle*{.5}}
\put(15,6){\circle*{.5}}
\put(15,-6){\circle*{.5}}
\end{picture}
\end{center}
\vspace{0.3cm}
\caption{Two possible strategies inside an $AB$-block: The path can either
move straight across or move along the interface for awhile and then move
across. Both strategies correspond to a coarse-grained step diagonally upwards.}
\label{fig-twostrat}
\end{figure}

\noindent
There is an explicit formula for $\kappa(a-c,1-b)$, which we will not specify
here. The idea behind Proposition~\ref{p:linkinf} is that the polymer follows the
$AB$-interface over a distance $bL$ during $cL$ steps and then wanders
away from the $AB$-interface to the diagonally opposite corner over a
distance $(1-b)L$ during $(a-c)L$ steps. The optimal strategy is obtained
by maximizing over $b$ and $c$ (see Fig.~\ref{fig-twostrat}). A similar
variational expression holds for $\psi_{BA}$.

With (\ref{psiAABB}) and (\ref{psiinflink}) we have identified the four block
pair free energies in terms of the \emph{single linear interface free energy}
$\phi^\cI$. This constitutes a \emph{major simplification}, in view of the methods
and techniques that are available for linear interfaces. We refer the reader to
the recent monograph by Giacomin~\cite{Gi07}, which describes a body of
mathematical ideas, techniques and results for copolymers in the vicinity
of linear interfaces.

\subsection{Percolation set}
\label{S1.4}

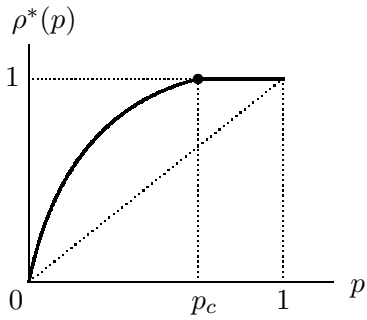
\begin{figure}
\begin{center}
\setlength{\unitlength}{0.45cm}
\begin{picture}(8,8)(0,1)
\put(0,0){\line(9,0){9}}
\put(0,0){\line(0,7){7}}
{\thicklines
\qbezier(5,6)(6.25,6)(7.5,6)
\qbezier(0,0)(1,5)(5,6)
}
\qbezier[40](5,0)(5,3)(5,6)
\qbezier[40](0,6)(2.5,6)(5,6)
\qbezier[40](7.5,0)(7.5,3)(7.5,6)
\qbezier[60](0,0)(3.75,3)(7.5,6)
\put(-.6,-.8){$0$}
\put(4.8,-.8){$p_c$}
\put(5,6){\circle*{.35}}
\put(9.5,-.3){$p$}
\put(-.5,7.5){$\rho^*(p)$}
\put(-.7,5.8){$1$}
\put(7.3,-.8){$1$}
\end{picture}
\vspace{1cm}
\end{center}
\caption{Qualitative picture of $p\mapsto\rho^*(p)$.}
\label{fig-rho*def}
\end{figure}

Let
\begin{equation}
\label{rho*def}
\rho^*(p) = \sup_{(\rho_{kl})\in\cR(p)} [\,\rho_{AA}+\rho_{AB}\,].
\end{equation}
This is the maximal frequency of $A$-blocks crossed by an infinite coarse-grained
path. The graph of $p\mapsto\rho^*(p)$ is sketched in Fig.~\ref{fig-rho*def}. For
$p \geq p_c$ the oil blocks percolate, and the maximal time spent in the oil by
a coarse-grained path is $1$. For $p<p_c$, on the other hand, the oil blocks do
not percolate and the maximal time spent in the oil is $<1$. For $(\alpha,\beta)
\in\CONE$, the copolymer prefers the oil over the water. Hence, the behavior
of the copolymer changes at $p=p_c$.

\section{Phase diagram for $p \geq p_c$}
\label{S2}

\begin{figure}
\begin{center}
\includegraphics[scale = 0.4]{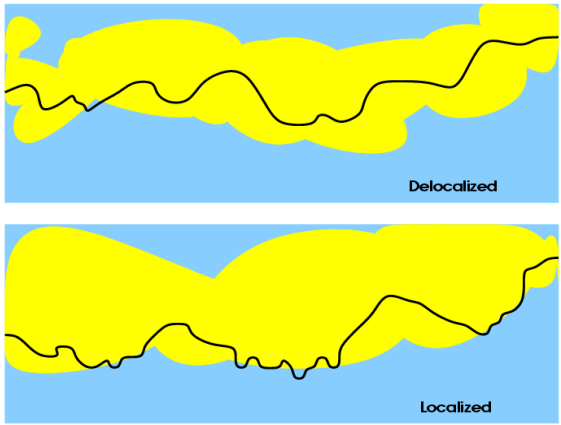}
\end{center}
\caption{Two possible strategies when the oil percolates.}
\label{fig-copolemulinfcl}
\end{figure}

The phase diagram is relatively simple in the supercritical regime. This is
because \emph{the oil blocks percolate}, and so the coarse-grained path can
choose between moving into the oil or running along the interface between
the oil and the water (see Fig.~\ref{fig-copolemulinfcl}).

The key result identifying the critical curve in the supercritical regime is
the following. Note that the criterion in (\ref{phinfcr}) is in terms of the
free energy of the single interface, and does not (!) depend on $p$.

\begin{proposition}
\label{p:phtrinfchar}
Let $p\geq p_c$. Then $(\alpha,\beta)\in\cL$ if and only if
\begin{equation}
\label{phinfcr}
\sup_{\mu \geq 1} \mu\left[\phi^{\cI}(\alpha,\beta;\mu)-\tfrac12\alpha
-\tfrac12\log 5\right] > \tfrac12\log\tfrac95.
\end{equation}
\end{proposition}

\noindent
Proposition~\ref{p:phtrinfchar} says that localization occurs if and only if the free
energy per step for the single linear interface exceeds the free energy per step for
an $AA$-block by a certain positive amount. This excess is needed to \emph{compensate}
for the loss of entropy that occurs when the path runs along the interface for awhile
before moving upwards from the interface to end at the diagonally opposite corner
(recall Fig.~\ref{fig-twostrat}). The constants $\frac12\log 5$ and $\frac12\log\frac95$
are special to our model. For the proof of Proposition~\ref{p:phtrinfchar} we refer
the reader to den Hollander and Whittington~\cite{dHoWh06}.

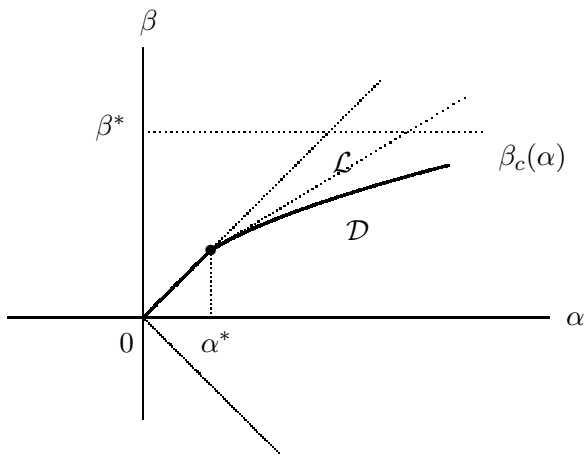
\begin{figure}
\begin{center}
\setlength{\unitlength}{0.45cm}
\begin{picture}(12,12)(0,-2)
\put(0,0){\line(12,0){12}}
\put(0,0){\line(0,8){8}}
\put(0,0){\line(0,-3){3}}
\put(0,0){\line(-4,0){4}}
{\thicklines
\qbezier(2,2)(4,3.2)(9,4.5)
\thicklines
\qbezier(0,0)(1,1)(2,2)
}
\qbezier[60](2,2)(4.5,3.5)(9.5,6.5)
\qbezier[60](0,5.5)(5,5.5)(10,5.5)
\qbezier[60](2,2)(4.5,4.5)(7,7)
\qbezier[15](2,0)(2,1)(2,2)
\qbezier[70](4,-4)(2,-2)(0,0)
\put(-.7,-1){$0$}
\put(12.5,-0.2){$\alpha$}
\put(-0.1,8.5){$\beta$}
\put(2,2){\circle*{.3}}
\put(1.7,-1){$\alpha^*$}
\put(-1.4,5.4){$\beta^*$}
\put(10.5,4.5){$\beta_c(\alpha)$}
\put(5.6,4.3){$\cal L$}
\put(6,2.3){$\cal D$}
\end{picture}
\end{center}
\vspace{1cm}
\caption{Qualitative picture of the phase diagram for $p \geq p_c$. There are
two phases, separated by a single critical curve $\alpha\mapsto\beta_c(\alpha)$.}
\label{fig-supcrit}
\end{figure}

With the help of Proposition~\ref{p:phtrinfchar} we can identify the supercritical
phase diagram. This runs via an analysis of the single interface free energy
$\phi^\cI$, for which we again refer to den Hollander and Whittington~\cite{dHoWh06}.
The phase diagram is sketched in Fig.~\ref{fig-supcrit}. The two phases are
characterized by
\begin{equation}
\begin{aligned}
\cD &= \{(\alpha,\beta)\in\CONE\colon\,
f(\alpha,\beta;p)=\tfrac12\alpha+\tfrac12\log 5\},\\
\cL &= \{(\alpha,\beta)\in\CONE\colon\,
f(\alpha,\beta;p)>\tfrac12\alpha+\tfrac12\log 5\},
\end{aligned}
\end{equation}
and are separated by a single critical curve $\alpha\mapsto\beta_c(\alpha)$.

The intuition behind the phase diagram is as follows. Pick a point $(\alpha,\beta)$
inside $\cD$. Then the polymer spends almost all of its time deep inside the
$A$-blocks. Increase $\beta$ while keeping $\alpha$ fixed. Then there will be
a larger energetic advantage for the polymer to move some of its monomers
from the $A$-blocks to the $B$-blocks by \emph{crossing the interface inside
the $AB$-block pairs}. There is some entropy loss associated with doing so,
but if $\beta$ is large enough, then the energy advantage will dominate,
so that $AB$-localization sets in. The value at which this happens depends
on $\alpha$ and is strictly positive. Since the entropy loss is finite,
for $\alpha$ large enough the energy-entropy competition plays out not only below
the diagonal, but also below a horizontal asymptote. On the other hand, for $\alpha$
small enough the loss of entropy dominates the energetic advantage, which is why
the critical curve has a piece that lies on the diagonal. At the critical value
$\alpha^*$, the critical curve has a slope discontinuity, because the linear
interface free energy is already strictly inside its localized region.. The
larger the value of $\alpha$ the larger the value of $\beta$ where $AB$-localization
sets in. This explains why the critical curve moves to the right and up.

In den Hollander and P\'etr\'elis~\cite{dHoPe07a} the following theorem is proved,
which completes the analysis of the phase diagram in Fig.~\ref{fig-supcrit}.

\begin{theorem}
\label{phtrsup}
Let $p \geq p_c$.\\
(i) $\alpha\mapsto\beta_c(\alpha)$ is strictly increasing on $[0,\infty)$.\\
(ii) For every $\alpha \in (\alpha^*,\infty)$ there exist $0<C_1<C_2<\infty$
and $\delta_0>0$ (depending on $p$ and $\alpha$) such that
\begin{equation}
\label{f2ndbds}
C_1\,\delta^2
\leq f\left(\alpha,\beta_c(\alpha)+\delta;p\right)
- f\left(\alpha,\beta_c(\alpha);p\right)
\leq C_2\,\delta^2 \qquad \forall\,\delta \in (0,\delta_0].
\end{equation}
(iii) $(\alpha,\beta)\mapsto f(\alpha,\beta;p)$ is infinitely differentiable
throughout $\cL$.
\end{theorem}

\noindent
Theorem~\ref{phtrsup}(i) implies that the critical curve never reaches the horizontal
asymptote, which in turn implies that $\alpha^*<\beta^*$ and that the slope at
$\alpha^*$ is $>0$. Theorem~\ref{phtrsup}(ii) shows that the phase transition
along the critical curve in Fig.~\ref{fig-supcrit} is \emph{second order off the
diagonal}. In contrast, we know that the phase transition is \emph{first order on
the diagonal}. Indeed, the free energy equals $\frac12\alpha+\frac12\log 5$ on and
below the diagonal segment between $(0,0)$ and $(\alpha^*,\alpha^*)$, and equals
$\frac12\beta+\frac12\log 5$ on and above this segment as is evident from interchanging
$\alpha$ and $\beta$. Theorem~\ref{phtrsup}(iii) tells us that the critical curve
in Fig.~\ref{fig-supcrit} is the only location in $\CONE$ where a phase transition
of finite order occurs.


\section{Phase diagram for $p<p_c$}
\label{S.CopolRa.4}

In the subcritical regime the phase diagram is \emph{much more complex} than
in the supcritical regime. The reason is that \emph{the oil does not percolate},
and so the copolymer no longer has the option of moving into the oil nor of
running along the interface between the oil and the water (in case it prefers
to localize). Instead, it has to every now and then cross blocks of water, even
though it prefers the oil.

\begin{figure}
\begin{center}
\setlength{\unitlength}{0.5cm}
\begin{picture}(12,12)(0,-3.5)
\put(0,0){\line(12,0){12}}
\put(0,0){\line(0,8){8}}
\put(0,0){\line(0,-3){3}}
\put(0,0){\line(-3,0){3}}
{\thicklines
\qbezier(0,0)(2,2)(4,4)
\qbezier(4,4)(3.5,2)(3,1)
\qbezier(3,1)(2,0)(1,-1)
\qbezier(3,1)(5,2)(9,3)
\qbezier(4,4)(6,5)(10,6)
}
\qbezier[20](4,0)(4,2)(4,4)
\qbezier[60](4,4)(6,6)(8,8)
\qbezier[40](0,0)(1.6,-1.5)(3.2,-3.0)
\qbezier[80](0,6.5)(4.5,6.5)(9,6.5)
\put(-.8,-.8){$0$}
\put(12.5,-0.2){$\alpha$}
\put(-0.1,8.5){$\beta$}
\put(3.7,-.8){$\alpha^*$}
\put(4,4){\circle*{.25}}
\put(3,1){\circle*{.25}}
\put(1.6,.8){$\cD_1$}
\put(6.5,1){$\cD_2$}
\put(6.5,3.5){$\cL_1$}
\put(8.5,7){$\cL_2$}
\end{picture}
\end{center}
\caption{Qualitative picture of the phase diagram for $p < p_c$. There are
four phases, separated by three critical curves, meeting at two tricritical
points.}
\label{fig-subcrit}
\end{figure}
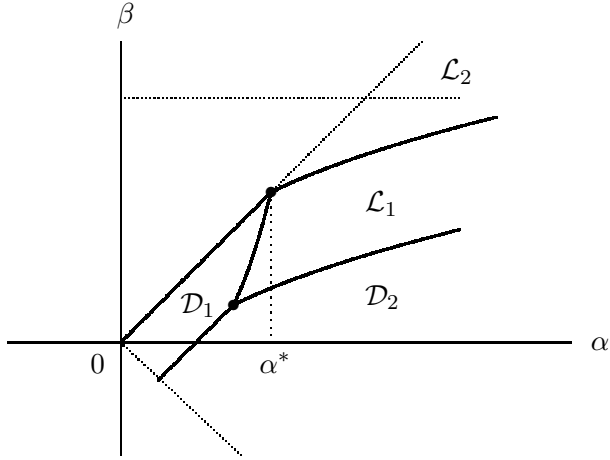

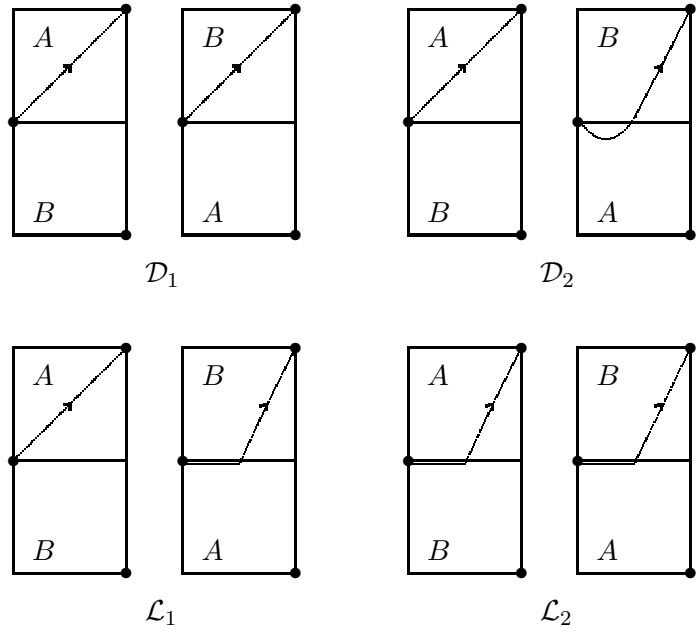
\begin{figure}
\begin{center}

\setlength{\unitlength}{0.25cm}
\begin{picture}(10,31)(12,-26)


{\thicklines
\qbezier(0,6)(3,6)(6,6)
\qbezier(6,0)(6,3)(6,6)
\qbezier(0,0)(3,0)(6,0)
\qbezier(6,-6)(6,-3)(6,0)
\qbezier(0,0)(0,3)(0,6)
\qbezier(0,0)(0,-3)(0,-6)
\qbezier(0,-6)(3,-6)(6,-6)
}
{\thicklines
\qbezier(2.7,3)(2.85,3)(3,3)
\qbezier(3,2.7)(3,2.85)(3,3)
}
\qbezier[60](0,0)(3,3)(6,6)
\put(1,4){$A$}
\put(1,-5.3){$B$}
\put(7,-8.5){$\cD_1$}

{\thicklines
\qbezier(9,6)(12,6)(15,6)
\qbezier(15,0)(15,3)(15,6)
\qbezier(9,0)(12,0)(15,0)
\qbezier(15,-6)(15,-3)(15,0)
\qbezier(9,0)(9,3)(9,6)
\qbezier(9,0)(9,-3)(9,-6)
\qbezier(9,-6)(12,-6)(15,-6)
}
{\thicklines
\qbezier(11.7,3)(11.85,3)(12,3)
\qbezier(12,2.7)(12,2.85)(12,3)
}
\qbezier[60](9,0)(12,3)(15,6)
\put(10,4){$B$}
\put(10,-5.3){$A$}

\put(0,0){\circle*{.5}}
\put(6,6){\circle*{.5}}
\put(6,-6){\circle*{.5}}
\put(9,0){\circle*{.5}}
\put(15,6){\circle*{.5}}
\put(15,-6){\circle*{.5}}


{\thicklines
\qbezier(21,6)(24,6)(27,6)
\qbezier(27,0)(27,3)(27,6)
\qbezier(21,0)(24,0)(27,0)
\qbezier(27,-6)(27,-3)(27,0)
\qbezier(21,0)(21,3)(21,6)
\qbezier(21,0)(21,-3)(21,-6)
\qbezier(21,-6)(24,-6)(27,-6)
}
{\thicklines
\qbezier(23.7,3)(23.85,3)(24,3)
\qbezier(24,2.7)(24,2.85)(24,3)
}
\qbezier[60](21,0)(24,3)(27,6)
\put(22,4){$A$}
\put(22,-5.3){$B$}
\put(28,-8.5){$\cD_2$}

{\thicklines
\qbezier(30,6)(33,6)(36,6)
\qbezier(36,0)(36,3)(36,6)
\qbezier(30,0)(33,0)(36,0)
\qbezier(36,-6)(36,-3)(36,0)
\qbezier(30,0)(30,3)(30,6)
\qbezier(30,0)(30,-3)(30,-6)
\qbezier(30,-6)(33,-6)(36,-6)
}
{\thicklines
\qbezier(34.15,2.9)(34.3,2.95)(34.45,3)
\qbezier(34.45,2.7)(34.45,2.85)(34.45,3)
}
\qbezier[50](33,0.2)(34.5,3.1)(36,6)
\qbezier[30](30,.2)(31.5,-2)(33,.2)
\put(31,4){$B$}
\put(31,-5.3){$A$}

\put(21,0){\circle*{.5}}
\put(27,6){\circle*{.5}}
\put(27,-6){\circle*{.5}}
\put(30,0){\circle*{.5}}
\put(36,6){\circle*{.5}}
\put(36,-6){\circle*{.5}}


{\thicklines
\qbezier(0,-12)(3,-12)(6,-12)
\qbezier(6,-18)(6,-15)(6,-12)
\qbezier(0,-18)(3,-18)(6,-18)
\qbezier(6,-24)(6,-21)(6,-18)
\qbezier(0,-18)(0,-15)(0,-12)
\qbezier(0,-18)(0,-21)(0,-24)
\qbezier(0,-24)(3,-24)(6,-24)
}
{\thicklines
\qbezier(2.7,-15)(2.85,-15)(3,-15)
\qbezier(3,-15.3)(3,-14.85)(3,-15)
}
\qbezier[60](0,-18)(3,-15)(6,-12)
\put(1,-14){$A$}
\put(1,-23.3){$B$}
\put(7,-26.5){$\cL_1$}

{\thicklines
\qbezier(9,-12)(12,-12)(15,-12)
\qbezier(15,-18)(15,-15)(15,-12)
\qbezier(9,-18)(12,-18)(15,-18)
\qbezier(15,-24)(15,-21)(15,-18)
\qbezier(9,-18)(9,-18)(9,-12)
\qbezier(9,-18)(9,-21)(9,-24)
\qbezier(9,-24)(12,-24)(15,-24)
}
{\thicklines
\qbezier(13.15,-15.1)(13.3,-15.05)(13.45,-15)
\qbezier(13.45,-15.3)(13.45,-15.15)(13.45,-15)
}
\qbezier[50](12,-18.2)(13.5,-14.9)(15,-12)
\qbezier[30](9,-18.2)(10.5,-18.2)(12,-18.2)
\put(10,-14){$B$}
\put(10,-23.3){$A$}

\put(0,-18){\circle*{.5}}
\put(6,-12){\circle*{.5}}
\put(6,-24){\circle*{.5}}
\put(9,-18){\circle*{.5}}
\put(15,-12){\circle*{.5}}
\put(15,-24){\circle*{.5}}


{\thicklines
\qbezier(21,-12)(24,-12)(27,-12)
\qbezier(27,-18)(27,-15)(27,-12)
\qbezier(21,-18)(24,-18)(27,-18)
\qbezier(27,-24)(27,-21)(27,-18)
\qbezier(21,-18)(21,-15)(21,-12)
\qbezier(21,-18)(21,-21)(21,-24)
\qbezier(21,-24)(24,-24)(27,-24)
}
{\thicklines
\qbezier(25.15,-15.1)(25.3,-15.05)(25.45,-15)
\qbezier(25.45,-15.3)(25.45,-15.15)(25.45,-15)
}
\qbezier[50](24,-18.2)(25.5,-14.9)(27,-12)
\qbezier[30](21,-18.2)(22.5,-18.2)(24,-18.2)
\put(22,-14){$A$}
\put(22,-23.3){$B$}
\put(28,-26.5){$\cL_2$}

{\thicklines
\qbezier(30,-12)(33,-12)(36,-12)
\qbezier(36,-18)(36,-15)(36,-12)
\qbezier(30,-18)(33,-18)(36,-18)
\qbezier(30,-18)(30,-15)(30,-12)
\qbezier(30,-18)(30,-21)(30,-24)
\qbezier(36,-18)(36,-21)(36,-24)
\qbezier(30,-24)(33,-24)(36,-24)
}
{\thicklines
\qbezier(34.15,-15.1)(34.3,-15.05)(34.45,-15)
\qbezier(34.45,-15.3)(34.45,-15.15)(34.45,-15)
}
\qbezier[50](33,-18.2)(34.5,-15.1)(36,-12)
\qbezier[30](30,-18.2)(31.5,-18.2)(33,-18.2)
\put(31,-14){$B$}
\put(31,-23.3){$A$}
\qbezier[50](33,-18.2)(34.5,-15.1)(36,-12)
\qbezier[30](30,-18.2)(31.5,-18.2)(33,-18.2)

\put(21,-18){\circle*{.5}}
\put(27,-12){\circle*{.5}}
\put(27,-24){\circle*{.5}}
\put(30,-18){\circle*{.5}}
\put(36,-12){\circle*{.5}}
\put(36,-24){\circle*{.5}}

\end{picture}

\end{center}
\caption{Behavior of the copolymer, inside the four block pairs containing
oil and water, for each of the four phases in Fig.~\ref{fig-subcrit}.}
\label{fig-subcritstrats}
\end{figure}

\begin{figure}
\begin{center}
 so\includegraphics[scale = 0.3]{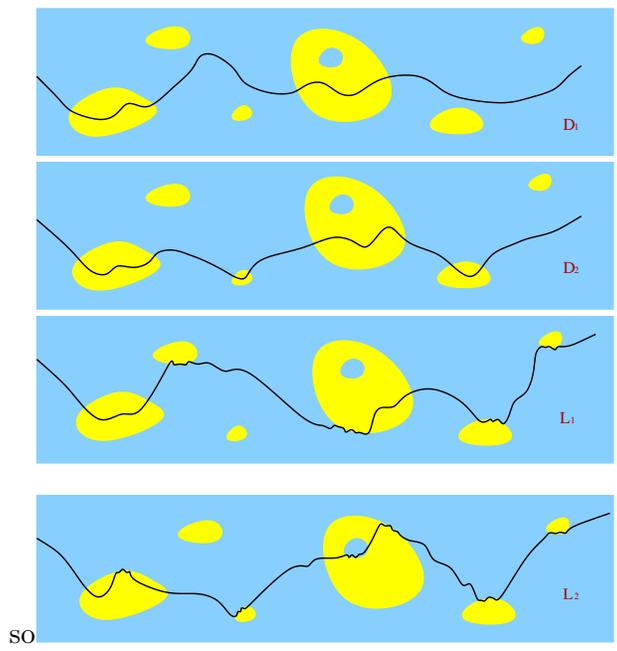}
\end{center}
\caption{Illustration of the four phases in Fig.~\ref{fig-subcrit} according
to Fig.~\ref{fig-subcritstrats}. Compare with Fig.~\ref{fig-copolemulinfcl}.}
\label{fig-copolemul4phas}
\end{figure}

It turns out that there are three (!) critical curves, all of which depend on
$p$. The phase diagam is sketched in Fig.~\ref{fig-subcrit}. For details on
the derivation, we refer to den Hollander and P\'etr\'elis~\cite{dHoPe07b}.
The copolymer has the following behavior in the four phases of Fig.~\ref{fig-subcrit},
as illustrated in Figs.~\ref{fig-subcritstrats} and \ref{fig-copolemul4phas}:
\begin{itemize}
\item[--]
$\cD_1$: \emph{fully delocalized} into $A$-blocks and $B$-blocks,
never inside a neighboring pair.
\item[--]
$\cD_2$: \emph{fully delocalized} into $A$-blocks and $B$-blocks,
sometimes inside a neighboring pair.
\item[--]
$\cL_1$: \emph{partially localized} near the interface in pairs of blocks
of which the $B$-block is crossed diagonally.
\item[--]
$\cL_2$: \emph{partially localized} near the interface in both types of
blocks.
\end{itemize}
(This is to be compared with the much simpler behavior in the two phases of
Fig.~\ref{fig-supcrit}, as given by Fig.~\ref{fig-twostrat}.)

The intuition behind the phase diagram is as follows. In $\cD_1$ and $\cD_2$,
$\beta$ is not large enough to induce localization. In both types of block
pairs, the reward for running along the interface is too small compared to
the loss of entropy that comes with having to cross the block at a steeper
angle. In $\cD_1$, where $\alpha$ and $\beta$ are both small, the copolmer
stays on one side of the interface in both types of block pairs. In $\cD_2$,
where $\alpha$ is larger, when the copolymer diagonally crosses a water block
(which it has to do every now and then because the oil does not percolate),
it dips into the oil block before doing the crossing. Since $\beta$ is small,
it still has no interest to localize. In $\cL_1$ and $\cL_2$, $\beta$ is large
enough to induce localization. In $\cL_1$, where $\beta$ is moderate, localization
occurs in those block pairs where the copolymer crosses the water rather
than the oil. This is because $\alpha>\beta$, making it more advantageous to
localize away from water than from oil. In $\cL_2$, where $\beta$ is large,
localization occurs in both types of block pairs.

Note that the piece between $\cD_1$ and $\cD_2$ is linear. This is because
in $\cD_1$ and $\cD_2$ the free energy is a function of $\alpha-\beta$ only.
The piece extends above the horizontal because no localization can occur
when $\beta\leq 0$. In $\cL_1$ and $\cL_2$ the free energy is a function of
$\alpha$ and $\beta$. Note that there are \emph{two tricritical points}, one
that depends on $p$ and one that does not.


Very little is known so far about the fine details of the four critical curves
in the subcritical regime. The reason is that in none of the four phases does
the free energy take on a simple form (contrary to what we saw in the supercritical
regime, where the free energy is simple in the delocalized phase). In particular, in
the subcritical regime there is no simple criterion like Proposition~\ref{p:phtrinfchar}
to characterize the phases. In den Hollander and P\'etr\'elis~\cite{dHoPe07b} it is
shown that the phase transition between $\cD_1$ and $\cD_2$ and between $\cD_1$ and
$\cL_1$ is second order, while the phase transition between $\cD_2$ and $\cL_1$ is
at least second order. It is further shown that the free energy is infinitely
differentiable in the interior of $\cD_1$ and $\cD_2$. The same is believed to be
true for $\cL_1$ and $\cL_2$, but a proof is missing.

It was argued in den Hollander and Whittington~\cite{dHoWh06} that \emph{the phase
diagram is discontinuous at $p=p_c$}. Indeed, none of the three critical curves
in the subcritical phase diagram in Fig.~\ref{fig-subcrit} converges to the
critical curve in the supercritical phase diagram in Fig.~\ref{fig-supcrit}.
This is because \emph{percolation versus non-percolation of the oil completely
changes the character of the phase transition(s)}.



\end{document}